# Bilateral binomial duplication formula

Martin Erik Horn, University of Potsdam
Am Neuen Palais 10, D - 14469 Potsdam, Germany
E-Mail: marhorn@uni-potsdam.de

**Abstract**

The bilateral binomial theorem with step width two gives a bilateral hypergeometric formula for $_2H_2\begin{bmatrix} a,\ a+1/2;\\ c,\ c+1/2;\end{bmatrix} z \end{bmatrix}$.

**Introduction**

At present nobody knows how the bilateral hypergeometric series

$$_2H_2\begin{bmatrix} a,\ b;\\ c,\ d;\end{bmatrix} z \end{bmatrix} = ? \qquad (1)$$

can be constructed. Only the Generalized Gauss Theorem [1, p. 181-183]

$$_2H_2\begin{bmatrix} a,\ b;\\ c,\ d;\end{bmatrix} 1 \end{bmatrix} = \Gamma\begin{bmatrix} c, d, 1-a, 1-b, c+d-a-b-1 \\ c-a, d-a, c-b, d-b \end{bmatrix}, \qquad (2)$$

1907 found by Dougall [2], [3, p.110] is known. A duplication of the parameters a and c of the bilateral binomial theorem [4]

$$_1H_1\begin{bmatrix} a;\\ c;\end{bmatrix} z \end{bmatrix} = \frac{(1-z)^{c-a-1}}{(-z)^{c-1}} \frac{\Gamma(c)\Gamma(1-a)}{\Gamma(c-a)} \qquad (3)$$

may give some hints about possible structural features of formula (1).

**Bilateral binomial theorem**

If integers of distance 1 in the Pascal Plane

$$\binom{x}{y} = \lim_{h \to 0} \frac{\Gamma(x+1+h)}{\Gamma(y+1+h)\Gamma(x-y+1+h)} \qquad (4)$$

are added, the bilateral binomial theorem

$$(1+z)^n = \sum_{k=-\infty}^{\infty} \binom{n}{k} z^k = \frac{\Gamma(n+1) \cdot z^k}{\Gamma(k+1)\Gamma(n-k+1)} {}_1H_1\begin{bmatrix} k-n;\\ k+1;\end{bmatrix} -z \end{bmatrix} \qquad (5)$$

resp.

$$(1-z)^n = \sum_{k=-\infty}^{\infty} \binom{n}{k} z^k = \frac{\Gamma(n+1) \cdot (-z)^k}{\Gamma(k+1)\Gamma(n-k+1)} {}_1H_1\begin{bmatrix} k-n;\\ k+1;\end{bmatrix} z \end{bmatrix} \qquad (6)$$

is obtained for $|z| = 1$ [5].

## First duplication of step width

The step width can be increased from distance one to distance two by adding formulae (5) and (6) according to

$$\frac{(1+z)^n}{z^K}+\frac{(1-z)^n}{(-z)^K} = \sum_{k=-\infty}^{\infty}\binom{n}{k}\left(\frac{z^k}{z^K}+\frac{(-z)^k}{(-z)^K}\right) \quad (7)$$

$$= \ldots +\binom{n}{k-2}\left(\frac{z^{k-2}}{z^K}+\frac{(-z)^{k-2}}{(-z)^K}\right)+\binom{n}{k-1}\left(\frac{z^{k-1}}{z^K}+\frac{(-z)^{k-1}}{(-z)^K}\right)$$

$$+\binom{n}{k}\left(\frac{z^k}{z^K}+\frac{(-z)^k}{(-z)^K}\right)+\binom{n}{k+1}\left(\frac{z^{k+1}}{z^K}+\frac{(-z)^{k+1}}{(-z)^K}\right)+\binom{n}{k+2}\left(\frac{z^{k+2}}{z^K}+\frac{(-z)^{k+2}}{(-z)^K}\right)+\ldots$$

Please notice that $k \in \mathbb{R}$ and that the left hand side of (7) does not depend on $k$. The parameter $K$ can be chosen at will.

With $K = k$

$$\frac{(1+z)^n}{z^k}+\frac{(1-z)^n}{(-z)^k} = \ldots + 2\binom{n}{k-2}z^{-2}+2\binom{n}{k}+2\binom{n}{k+2}z^2+2\binom{n}{k+4}z^4+\ldots \quad (8)$$

odd terms cancel out. Unfortunately the left hand side now depends on $k$. This formula can be transformed into the standard hypergeometric notation.

$$\frac{(1+z)^n}{z^k}+\frac{(1-z)^n}{(-z)^k}$$
$$= \frac{2\cdot\Gamma(n+1)}{\Gamma(k+1)\Gamma(n-k+1)}\left(\ldots+\frac{k(k-1)}{(n-k+1)(n-k+2)}z^{-2}+1+\frac{(n-k)(n-k-1)}{(k+1)(k+2)}z^2+\ldots\right) \quad (9)$$

$$= \frac{2\cdot\Gamma(n+1)}{\Gamma(k+1)\Gamma(n-k+1)}\left(\ldots+\frac{\frac{k}{2}\frac{k-1}{2}}{\frac{k-n-1}{2}\frac{k-n-2}{2}}z^{-2}+1+\frac{\frac{k-n}{2}\frac{k-n+1}{2}}{\frac{k+1}{2}\frac{k+2}{2}}z^2+\ldots\right) \quad (10)$$

$$= \frac{2\cdot\Gamma(n+1)}{\Gamma(k+1)\Gamma(n-k+1)}\,_2\mathrm{H}_2\!\left[\begin{array}{c}\frac{k-n}{2},\ \frac{k-n+1}{2};\\ \frac{k+1}{2},\ \frac{k+2}{2};\end{array} z^2\right] \quad (11)$$

Finally $k = 2c - 1$, $n = 2c - 2a - 1$ and $z^2 \to z$ give the bilateral binomial duplication formula

$$_2\mathrm{H}_2\!\left[\begin{array}{c}a,\ a+1/2;\\ c,\ c+1/2;\end{array} z\right] = \frac{\Gamma(2c)\Gamma(1-2a)}{2\cdot\Gamma(2c-2a)}\left(\frac{(1+\sqrt{z})^{2c-2a-1}}{\sqrt{z}^{2c-1}}+\frac{(1-\sqrt{z})^{2c-2a-1}}{(-\sqrt{z})^{2c-1}}\right) \quad (12)$$

This formula (12) shows a nice symmetry between $(\sqrt{z})$ and $(-\sqrt{z})$.

## Second duplication of step width

Alternatively the step width can be increased from distance one to distance two by subtracting formulae (5) and (6) according to

$$\frac{(1+z)^n}{z^K} - \frac{(1-z)^n}{(-z)^K} = \sum_{k=-\infty}^{\infty} \binom{n}{k}\left(\frac{z^k}{z^K} - \frac{(-z)^k}{(-z)^K}\right) \quad (13)$$

$$= ... + \binom{n}{k-2}\left(\frac{z^{k-2}}{z^K} - \frac{(-z)^{k-2}}{(-z)^K}\right) + \binom{n}{k-1}\left(\frac{z^{k-1}}{z^K} - \frac{(-z)^{k-1}}{(-z)^K}\right)$$

$$+ \binom{n}{k}\left(\frac{z^k}{z^K} - \frac{(-z)^k}{(-z)^K}\right) + \binom{n}{k+1}\left(\frac{z^{k+1}}{z^K} - \frac{(-z)^{k+1}}{(-z)^K}\right) + \binom{n}{k+2}\left(\frac{z^{k+2}}{z^K} - \frac{(-z)^{k+2}}{(-z)^K}\right) + ...$$

Please notice that $k \in \mathbb{R}$ and that the left hand side of (13) does not depend on $k$. The parameter $K$ can be chosen at will.

With $K = k$

$$\frac{(1+z)^n}{z^k} - \frac{(1-z)^n}{(-z)^k} = ... + 2\binom{n}{k-3}z^{-3} + 2\binom{n}{k-1}z^{-1} + 2\binom{n}{k+1}z^1 + 2\binom{n}{k+3}z^3 + ... \quad (14)$$

even terms cancel. Unfortunately the left hand side now depends on $k$. This formula can be transformed into the standard hypergeometric notation.

$$\frac{(1+z)^n}{z^k} - \frac{(1-z)^n}{(-z)^k}$$

$$= \frac{2\cdot\Gamma(n+1)}{\Gamma(k+1)\Gamma(n-k+1)}\left(... + \frac{k(k-1)(k-2)}{(n-k+1)(n-k+2)(n-k+3)}z^{-3} + \frac{k}{n-k+1}z^{-1}\right. \quad (15)$$

$$\left. + \frac{n-k}{k+1}z^1 + \frac{(n-k)(n-k-1)(n-k-2)}{(k+1)(k+2)(k+3)}z^3 + ...\right)$$

$$= \frac{2}{z}\frac{\Gamma(n+1)}{\Gamma(k)\Gamma(n-k+2)}\left(... + \frac{(k-1)(k-2)}{(n-k+2)(n-k+3)}z^{-2} + 1 + \frac{(n-k)(n-k-1)}{k(k+1)}z^2\right.$$

$$\left. + \frac{(n-k)(n-k-1)(n-k-2)(n-k-3)}{k(k+1)(k+2)(k+3)}z^4 + ...\right) \quad (16)$$

$$= \frac{2}{z}\frac{\Gamma(n+1)}{\Gamma(k)\Gamma(n-k+2)}\left(... + \frac{\frac{k-1}{2}\frac{k-2}{2}}{\frac{k-n-2}{2}\frac{k-n-3}{2}}z^{-2} + 1 + \frac{\frac{k-n-1}{2}\frac{k-n}{2}}{\frac{k}{2}\frac{k+1}{2}}z^2 + ...\right) \quad (17)$$

$$= \frac{2}{z}\frac{\Gamma(n+1)}{\Gamma(k)\Gamma(n-k+2)} {}_2H_2\left[\begin{array}{c}\frac{k-n-1}{2}, \frac{k-n}{2};\\ \frac{k}{2}, \frac{k+1}{2};\end{array} z^2\right] \quad (18)$$

Finally $k = 2c$, $n = 2c - 2a - 1$ and $z^2 \to z$ give again the bilateral binomial duplication formula

$$_2H_2\begin{bmatrix} a,\ a+1/2; \\ c,\ c+1/2; \end{bmatrix} z \end{bmatrix} = \frac{\Gamma(2c)\Gamma(1-2a)}{2\cdot\Gamma(2c-2a)}\left(\frac{(1+\sqrt{z})^{2c-2a-1}}{\sqrt{z}^{2c-1}} + \frac{(1-\sqrt{z})^{2c-2a-1}}{(-\sqrt{z})^{2c-1}}\right), \quad (12)$$

and again it does not matter whether the square root of $z$ is chosen positive or negative.

**Special cases**

As expected formulae (2) and (12) coincide at

$$_2H_2\begin{bmatrix} a,\ a+1/2; \\ c,\ c+1/2; \end{bmatrix} 1 \end{bmatrix} = \frac{\Gamma(2c)\Gamma(1-2a)}{\Gamma(2c-2a)}\cdot 2^{2c-2a-2} \quad (19)$$

which can be shown by using Legendre's duplication formula [3, p.22].

For $z = -1$ the bilateral binomial duplication formula (12) gives

$$_2H_2\begin{bmatrix} a,\ a+1/2; \\ c,\ c+1/2; \end{bmatrix} -1 \end{bmatrix} = \frac{\Gamma(2c)\Gamma(1-2a)}{2\cdot\Gamma(2c-2a)}\left(\frac{(1+i)^{2c-2a-1}}{i^{2c-1}} + \frac{(1-i)^{2c-2a-1}}{(-i)^{2c-1}}\right) \quad (20)$$

$$= \frac{\Gamma(2c)\Gamma(1-2a)}{\Gamma(2c-2a)} 2^{c-a-\frac{3}{2}}\left(\left(\cos\frac{\pi}{4}+i\sin\frac{\pi}{4}\right)^{-2a}\left(\cos\frac{\pi}{4}-i\sin\frac{\pi}{4}\right)^{2c-1} \right.$$
$$\left. +\left(\cos\frac{\pi}{4}-i\sin\frac{\pi}{4}\right)^{-2a}\left(\cos\frac{\pi}{4}+i\sin\frac{\pi}{4}\right)^{2c-1}\right) \quad (21)$$

$$= \frac{\Gamma(2c)\Gamma(1-2a)}{\Gamma(2c-2a)} 2^{c-a-\frac{3}{2}}\left(e^{-i\frac{\pi}{2}\left(c+a-\frac{1}{2}\right)} + e^{i\frac{\pi}{2}\left(c+a-\frac{1}{2}\right)}\right) \quad (22)$$

$$= \frac{\Gamma(2c)\Gamma(1-2a)}{\Gamma(2c-2a)} 2^{c-a-\frac{3}{2}} \cos\left(\frac{2c+2a-1}{4}\pi\right) \quad (23)$$

For $c = 1/2$ formula (23) transforms into a unilateral hypergeometric series

$$_2H_2\begin{bmatrix} a,\ a+1/2; \\ 1/2,\ 1; \end{bmatrix} -1 \end{bmatrix} = {}_2F_1\begin{bmatrix} a,\ a+1/2; \\ 1/2; \end{bmatrix} -1 \end{bmatrix} = 2^{-a}\cos\frac{\pi a}{2} \quad (24)$$

and can be compared with Kummers hypergeometric sum [3, p.126]

$$_2F_1\begin{bmatrix} a,\ b; \\ a-b+1; \end{bmatrix} -1 \end{bmatrix} = \frac{\Gamma(a-b+1)\Gamma(b)}{\Gamma(a+1)\Gamma(a/2-b+1)} \quad (25)$$

for $b = a + 1/2$

$$_2F_1\begin{bmatrix} a,\ a+1/2; \\ 1/2; \end{bmatrix} -1 \end{bmatrix} = \frac{\Gamma(1/2)\Gamma(a/2+1)}{\Gamma(a+1)\Gamma(-a/2+1/2)} = 2^{-a}\cos\frac{\pi a}{2} \quad (26)$$

using Legendre's duplication formula [3, p.22].

**Outlook**

Using the same duplication procedure again and again bilateral hypergeometric series of higher rank can be constructed.